\input amstex
\documentstyle{amsppt}
\magnification=\magstep1
\baselineskip 12 true pt
\parskip=0pt plus 3pt
\hsize16truecm \hoffset0.4truecm \vsize24truecm \voffset0.2truecm
\raggedbottom
\parindent=0.5cm

\def\cd{\cdots}

\def\max{\text{\rm max}\, }

\def\dim{\text{\rm dim}\, }

\def\al{\alpha}

\def\la{\lambda}

\def\U{\frak U}

\def\g{\frak g}

\def \bZ{\Bbb Z}
\def \bN{\Bbb N}

\def \bC{\Bbb C}

\font\boldtitlefont=cmb10 scaled\magstep1

\font\bigmath=cmbxti10 scaled\magstep2
\font\smmath=cmr10

\font\scten=cmti10

\leftheadtext {Ye and Zhou}
\rightheadtext {A New Multiplicity Formula}

\topmatter
\title
{\boldtitlefont  A New Multiplicity Formula for the Weyl Modules
of Type ${\hbox{\bigmath A}}^\dag$}
\endtitle
\thanks
{\dag \; Supported in part by the NNSFC (10271088).}
\endthanks
\author
{Ye Jiachen$^{1}$~\;~Zhou Zhongguo$^{1,~2}$}
\endauthor
\affil
{\scten $^1$Department of Applied Mathematics, Tongji University,\\
Shanghai 200092, People's Republic of China\\
E-Mail: jcye\@mail.tongji.edu.cn\\
\medskip
$^2$Department of Mathematics, Jinan University,\\
Jinan  Shandong 250022, People's Republic of China\\
E-Mail: zhgzhou\@21cn.com\\}
\endaffil
\endtopmatter
\noindent {\bf Abstract:} \quad {\smmath A monomial basis and a
filtration of subalgebras for the universal enveloping algebra
$\U(\g)$ of a complex simple Lie algebra $\frak g_l$ of type
$A_l$ is given in this note. In particular, a new multiplicity
formula for the Weyl module $V(\la)$ of $\U(\g_l)$ is obtained in
this note.}
\medskip
\noindent {\bf 2000 MR Subject Classification:} 17B10 \quad 20G05
\vskip0.4cm
\document
Let $\g_l$ be a complex simple Lie algebra of type $A_l$, and
$\U=\U(\g_l)$ its universal enveloping algebra. For any dominant
integer weiget $\la\in \Lambda_+$, $V(\la)$ denotes a finite
dimensional irreducible $\U(\g_l)$-module, the Weyl module.
Following Littelmann [2], we define a new monomial basis and a
filtration of subalgebras for $\U(\frak g_l)$. Furthermore, we
obtain a new basis and a new multiplicity formula for the Weyl
module $V(\la)$ of $\U(\g_l)$ in this note.

This paper is organized as follows: we introduce an ordering
relation on $\U(\g_l)^-$ in the first section; we define a new
basis of $\U(\g_l)^-$ in Section 2; we record some useful
commutative formulas and construct a filtration of subalgebras of
$\U(\g_l)$ in Section 3; our main results concerning a new basis
and a new multiplicity formula for the Weyl module $V(\la)$ of
$\U(\g_l)$ is given in Section 4; several examples for $\g_l$
being of type $A_2$ and $A_4$ is given in the last section. We
shall freely use the notations in [1] without further comments.

We believe that our method could be generalized to the case of
$D_l$ at least. Moreover, our results may also be generalized to
the cases of $B_l$ and $C_l$. We will deal with them in a further
note.

\head {1. An ordering relation on $\U(\g_l)^-$}\endhead

{\bf 1.1.} Let $\Delta=\{\al_1, \al_2,\cd, \al_l\}$ be the set of
simple roots. Set $\al_{i~i}=\al_{i}$, and
$$\al_{i~j}=\al_i+\al_{i+1}+\cd+\al_j, \qquad \text {with}\qquad 1\le i\le j\le l.$$
Then
$$\Phi^+=\{\al_{i~j}, \quad 1\le i < j\le l\}$$ is the set of positive
roots which has $\dfrac{1}{2}l(l+1)$ elements. Fix an ordering of
positive roots as follows:
$$\al_1, \al_2, \al_{1~2},\cdots, \al_i, \al_{i-1~i},\cdots,
\al_{1~i},\cdots, \al_l,\al_{l-1~l},\cdots, \al_{1~l}.$$ Then
define a Chevalley basis of $\g_l$
$$e_i=e_{\al_i}, e_{i~j}=e_{\al_{i~j}}, f_i=f_{\al_i}, f_{i~j}=f_{\al_{i~j}},
h_i, \quad 1\le i \le l,\;\; 1\le i<j\le l,$$ accordingly. Let
$\bN$ be the set of non-negative integers. The Kostant $\bZ$-form
$\U_{\bZ}$ of $\U$ is the $\bZ$-subalgebra of $\U$ generated by
the elements $e_{\al}^{(k)}: =e_{\al}^k/k!,
f_{\al}^{(k)}:=f_{\al}^k/k!$ for $\al \in \Phi^+$ and $k \in \bN$.
Set
$$
{{h_i+c}\choose k}: = \frac{(h_i+c)(h_i+c-1)\cdots
(h_i+c-k+1)}{k!}.
$$
Then $ {{h_i+c}\choose k}\in \U_{\bZ}$ for $1\le i\le l, c\in\bZ,
k\in\bN$. Let $\U^+, \U^-, \U^0$ be the positive part, negative
part and zero part of $\U$, respectively. They are generated by
$e_{\al}^{(k)}$, $f_{\al}^{(k)}$ and ${{h_i} \choose k}$ with
$k\in\bN, \al\in \Phi^+$ and $1\le i \le l$, respectively. The
algebra $\U$ is a Hopf algebra which has a triangular
decomposition $\U=\U^-\U^0\U^+$. It is known that the PBW-type
basis for $\U$ has the form of
$$
\prod_{\al\in R_+}f_{\al}^{(a_\al)}\prod_{i=1}^l{{h_i}\choose b_i}
\prod_{\al\in R_+}e_{\al}^{(c_\al)}
$$
with $a_\al, b_i, c_\al \in \bN$. In particular, if we define
$$I=(i_1, i_2, \cdots, i_{\frac{l(l+1)}{2}})\in
\bN^{\frac{l(l+1)}{2}},$$ then
$$f^I=f_1^{(i_1)}f_2^{(i_2)}f_{1~2}^{(i_3)}\cdots f_{1~l}^{(i_{\frac{l(l+1)}{2}})}$$
forms a PBW-type basis of $\U^-$ with all $I\in \bN^{\frac
{l(l+1)}{2}}$. In particular, one has $f^{\text{\bf 0}}=1$ when
$I=(0, 0, \cdots, 0)=\text{\bf 0}$.
\par
{\bf 1.2.} First of all, we define an ordering on
$\bN^{\frac{l(l+1)}{2}}$ \lq\lq$\prec$\rq\rq as follows: for any
$I, I'\in \bN^{\frac{l(l+1)}{2}}$, $I=(i_1, i_2, \cdots,
i_{\frac{l(l+1)}{2}})$ and $I'=(i'_1, i'_2, \cdots,
i'_{\frac{l(l+1)}{2}})$, if there exists a $k$ with $1\le k \le
\frac{l(l+1)}{2}$ such that $i_k<i'_k$ and $i_j=i'_j$ for all
$j>k$, then we say $I \prec I'$; otherwise, one has $I = I'$. It
is easy to see that for any $I, I'\in \bN^{\frac{l(l+1)}{2}}$, if
$I\ne I'$, we must have either $I \prec I'$ or $I' \prec I$.
Therefore, we can define an ordering on $\U^-$ ~$\lq\lq \prec "$
in the same way: we say $f^I \prec f^{I'}$ if and only if $I \prec
I'$. It is easy to see that different basis elements do not be
equal. Any element in $\U^-$ can be written uniquely in terms of
$$f=\sum_{I\in \bN^{\frac{l(l+1)}{2}}} a_{{_I}} f^I, \qquad \text {with}
\qquad a_{{_I}} \in \bC.$$ Moreover, we can define the leading
element of $f$ $\max f=f^I$ when all the other $f^{I'}\prec f^I$
with $a_{_{I'}}\ne 0$. Therefore, one has the following claim:

{\bf 1.3.} If $f_1, f_2, \cdots $, $f_m \in \U^-$ with $\max
f_1\prec \max f_2\prec \cdots \prec \max f_m$. Then $f_1$, $f_2$,
$\cdots $, $f_m$ are linearly independent.

\head {2. Some commutator formulas and \\a class of special
subalgebras in $\U(\g_l)$}\endhead

{\bf 2.1.} For $1\le i, j\le l$, one has the following commutator
formulas (cf. [1]).
$$\align
&e_if_j=f_je_i, \quad {\text {\rm when}} \;\; i\ne j;\tag 1\\
&e_i^{(a)}f_i^{(b)}=\sum_{k=0}^{\text {min} (a,b)}
f_i^{(b-k)}{{h_i-a-b+2k}\choose k}e_i^{(a-k)};\tag 2\\
&h_if_j^{(k)}=f_j^{(k)}h_i-k\al_j(h_i)f_j^{(k)};\tag 3\\
&{{h_i+a}\choose b}f_{j}^{(k)}=f_j^{(k)}{{h_i-k\al_j(h_i)+a}
\choose b};\tag 4\\
&e_if_l^{(a_l)}\cdots f_i^{(a_i)} \cdots f_{1}^{(a_1)}
=f_l^{(a_l)}\cdots f_i^{(a_i)}\cdots f_{1}^{(a_1)}e_i+\tag 5\\
&\;\;+f_1^{(a_l)}\cdots f_i^{(a_i-1)}\left(h_i-a_i+1\right)
f_{i-1}^{(a_{i-1})}\cdots f_{1}^{(a_1)}=f_l^{(a_l)}\cdots
f_i^{(a_i)}\cdots f_{1}^{(a_1)}e_i +\\
&\;\;+f_1^{(a_l)}\cdots f_i^{(a_i-1)} \cdots f_{1}^{(a_1)}
\left(h_i-a_i+1-\sum_{k=1}^{i-1}a_k\al_k(h_i)\right).
\endalign$$

{\bf 2.2.} Furthermore, elements $f_i, f_{i~j},$ ($1 \le i < j \le
l$) satisfy the following commutator relations (cf. [4] or [5]):
$$ \align
&f_{i+1}f_{i}=f_{i}f_{i+1}+f_{i~i+1};\tag 1\\
&f_if_{j}=f_{j}f_i,\quad {\text {\rm when}} \;\; |i-j|\ne 1;\tag 2\\
&f_{i+1~j}f_i=f_if_{i+1~j}+f_{i~j}\;\;
{\text {\rm \ or}}\;\;f_{j+1}f_{i~j}=f_{i~j}f_{j+1}+f_{i~j+1};\tag 3\\
&f_{i~j}f_{k}=f_{k}f_{i~j},\;\;\;\;\;\; {\text {\rm when}} \;\;
i-k \ne 1\;\; {\text {\rm \ or}}\;\;k-j\ne 1;\tag 4\\
&f_{j+1~k}f_{i~j}=f_{i~j}f_{j+1~k}+f_{i~k};\tag 5\\
&f_{i~j}f_{k~h}=f_{k~h}f_{i~j}, \;\;\; {\text {\rm when}} \;\; k-j
\ne 1\;\; {\text {\rm \ or}}\;\; i-h \ne 1;\tag 6\\
&f_{i+1}^{(a)}f_i^{(b)}=\sum_{k=0}^{\text {min} (a,b)}
f_{i~i+1}^{(k)}f_i^{(b-k)}f_{i+1}^{(a-k)}, \quad 1\le
i\le l-1.\tag 7\\
\endalign $$

{\bf 2.3.} Let us construct a class of special subalgebras
$\U(\g_i), 1\le i\le l$, of $\U$ as follows. Set
$$\U(\g_i)=\langle e_j^{(a_j)}, f_j^{(b_j)}, {{h_j+c}\choose
k}|~a_j, b_j, c, k \in \bN, 1\le j \le i \rangle.$$ Then one has
$$0\subseteq \U(\g_1)\subseteq \U(\g_2)\subseteq \cdots \subseteq
\U(\g_l)=\U.$$ The set of positive roots in $\U(\g_i)$ is just
that of the first $\dfrac{1}{2}i(i+1)$ roots according to the
ordering of $\Phi^+$.

\head {3. A monomial basis of $\U(\g_l)^-$}\endhead

{\bf 3.1.} Let $K=(k_1^l, k_2^{l-1}, k_1^{l-1}, \cdots,
k_i^{l-i+1}, k_{i-1}^{l-i+1}, \cdots, k_1^{l-i+1}, \cdots, k_l^1,
k_{l-1}^1$, $\cdots, k_1^1)$ $\in \bN^{\frac{l(l+1)}{2}}$. Define
an index set
$$\Pi:=\{K \in \bN^{\frac{l(l+1)}{2}}|~k_i^{l-i+1}\ge
k_{i-1}^{l-i+1}\ge \cdots \ge k_1^{l-i+1}, \; 1\le i \le l\}.$$
For any $K\in \Pi$, one has such a monomial
$$\align
\theta^K=&f_1^{(k_1^l)}f_2^{(k_2^{l-1})}f_1^{(k_1^{l-1})}\cdots
f_i^{ (k_i^{l-i+1})}f_{i-1}^{(k_{i-1}^{l-i+1})}\cdots
f_1^{(k_1^{l-i+1})} \cdots \\
&f_l^{(k_l^1)}f_{l-1}^{(k_{l-1}^1)}\cdots f_1^{(k_1^1)}\in \U^-.
\endalign$$

The following theorem was first proved by Littelmann [2].

{\bf 3.2.} {\bf {Theorem}}\;\;{\it {The set $\{\theta^K |~K\in \Pi\}$
forms a basis of the $\bZ$-form of $\U^-$.}}

{\bf{Proof}\;\;} First of all, we have to show that elements of
the set $\{\theta^K |~K\in \Pi\}$ are linearly independent.

Since $\{f^I|~I\in \bN^{\frac{l(l+1)}{2}}\}$ forms a PBW-type
basis of $\U^-$, one has for any $K\in\Pi$, $\theta^K\in \U^-$ and
$$\theta^K=\sum_I a_{{_I}}f^I, \quad a_{{_I}}\in \bZ, \;\;I\in
\bN^{\frac{l(l+1)}{2}}.$$

Moreover, for any $K\in \Pi$, one has $I(K)=(k_1^l, k_2^{l-1}-
k_1^{l-1}, k_1^{l-1}, \cdots, k_i^{l-i+1}-k_{i-1}^{l-i+1},
k_{i-1}^{l-i+1}-k_{i-2}^{l-i+1}, \cdots,
k_{2}^{l-i+1}-k_{1}^{l-i+1}, k_1^{l-i+1}, \cdots,
k_l^{1}-k_{l-1}^{1}, k_{l-1}^{1}-k_{l-2}^{1},\cdots,
k_{2}^{1}-k_{1}^{1}, k_1^1)\in \bN^{\frac{l(l+1)}{2}}$, because
$k_i^j\ge k_{i-1}^j$ for all $1\le i, j \le l$ with $i+j\le l+1$.
It is easy to calculate that
$$\max \theta^K=f^{I(K)}, \quad \text {with
coefficient} \;\;1.$$ Therefore, one has
$$\theta^K=f^{I(K)}+\sum_{I\prec I(K)}a_{{_I}}f^I.$$

Note the fact that various $\theta^K$s and $\max \theta^K$s are
different, when the corresponding $K$s are different. We can
conclude that elements of the set $\{\theta^K |~K\in \Pi\}$ are
linearly independent.

Next we show that the set $\{\theta^K |~K\in \Pi\}$ generate
$\U^-_{\bZ}$. For any $I=(i_1^l, i_2^{l-1}, i_1^{l-1}$, $\cdots$,
$i_i^{l-i+1}, i_{i-1}^{l-i+1}, \cdots, i_1^{l-i+1}, \cdots$,
$i_l^1, i_{l-1}^1$, $\cdots, i_1^1) \in \bN^{\frac{l(l+1)}{2}}$,
we define $K(I)=(i_1^l, i_2^{l-1}+i_1^{l-1}, i_1^{l-1}, \cdots,
\sum_{p=1}^{j} i_p^{l-j+1}, \cdots, i_{2}^{l-j+1}+i_{1}^{l-j+1},
i_1^{l-j+1}, \cdots, \sum_{p=1}^l i_p^{1},
\sum_{p=1}^{l-1}i_{p}^{1}, \cdots, i_{2}^{1}+i_{1}^{1}, i_1^1) \in
\Pi.$ Then one has
$$\theta^{K(I)}=f^I+\sum_{I'\prec I}a_{_{I'}}f^{I'}.$$ An easy
induction on the ordering of $\bN^{\frac{l(l+1)}{2}}$ shows that
$$f^I=\theta^{I(K)}+\sum_{K'\in \Pi}c_{_{K'}}\theta^{K'}\quad \text
{with} \;\;c_{_{K'}}\in \bZ.$$

Combining the above facts, we show that the set $\{\theta^K |~K\in
\Pi\}$ forms a basis of the $\bZ$-form of $\U^-$.

{\bf 3.3.} Define $\Pi_{l-1}:=\{ K\in \Pi|~ k_j^1=0, \; 1\le j\le
l\}\subseteq \Pi$. We can see from the above discussion that the
set $\{\theta^K |~ K\in \Pi_{l-1}\}$ forms a basis of the
$\bZ$-form of $\U(\g_{l-1})^-$.

Set $\Pi':=\{K\in \Pi|~ k_j^i=0, \; 1<i\le l\}$.

{\bf 3.4.} If we define the ordinary vector addition in $\Pi$, one
has the following claims:
$$\align &\text {(1)} \;\;\Pi=\Pi_{l-1}\oplus \Pi';\\
&\text {(2)}\;\;\text {If}\; K_2\in \Pi_{l-1}\;\text {and} \;
K_1\in \Pi', \;\text {then}
\; \theta^{K_2}\theta^{K_1}=\theta^{K_2+K_1};\\
&\text {(3)}\;\;\text {If} \; K_1, K_1'\in \Pi'\; \text {with}\;
K_1\prec K_1',\; \text {then}\; K_2+K_1\prec K_1' \;\text
{for any}\; K_2\in \Pi_{l-1}.\\
\endalign$$

\head {4. A new multiplicity formula of $V(\la)$}\endhead

{\bf 4.1.} Let $\Lambda$ be the set of weights for $\g_l$, and
$\omega_1, \omega_2, \cdots, \omega_l$ the set of fundamental
dominant weights. Then the set of dominant weights $\Lambda^+$ is
defined to be
$$\{\la=(\la_1, \la_2, \cdots, \la_l)=\sum_{i=1}^l \la_i\omega_i \;\;
{\text {\rm with all}} \;\;\la_i\in \bN\}.$$

Let $E$ be the real vector space spanned by $\al_1, \al_2,\cd,
\al_l$. It is well-known that $\al_1^{\vee}, \al_2^{\vee},\cdots,
\al_l^{\vee}$ again form a basis of $E$, and $\omega_1, \omega_2,
\cdots, \omega_l$ form the dual basis relative to the inner
product on $E$: $(\omega_i,\al_j^{\vee})=\dfrac{2(\omega_i,
\al_j)}{(\al_j,\al_j)}=\delta_{i~j}$. If we restrict ourselves to
considering the $(l-1)$-dimensional subspaces $E'$ of $E$ spanned
by $\al_1, \al_2,\cd, \al_{l-1}$, then $\al_1^{\vee},
\al_2^{\vee}, \cdots, \al_{l-1}^{\vee}$ and $\omega_1, \omega_2,
\cdots, \omega_{l-1}$ remain the dual bases of $E'$ relative to
the inner product on $E$. Therefore, we can consider the
restriction of $\U(\g_{l})$ to $\U(\g_{l-1})$, and the restriction
of $\la =(\la_1, \la_2, \cdots, \la_l)$ as a weight of $\g_l$ to
$\la_{\g_{l-1}} =(\la_1, \la_2, \cdots, \la_{l-1})$ as a weight of
$\g_{l-1}$. Moreover, let $\la =(\la_1, \la_2, \cdots, \la_l) \in
\Lambda^+$ be a dominant weight, and $v$ a maximal vector of
weight $\la$ of the $\U(\g_l)$-module $V(\la)$. Then
$V(\la)|_{\U(\g_{l-1})}$ denotes the restriction of $V(\la)$ to a
$\U(\g_{l-1})$-module.

We can make use of the recursive property of the basis $\{\theta^K
|~ K\in \bN^{\frac{l(l+1)}{2}}\}$ to construct a new basis of the
finite-dimensional irreducible $\g$-module $V(\la)$ with $\la \in
\Lambda^+$, and to get a new multiplicity formula of $V(\la)$.
Following Littermann [2], we define $\la_i^j$ in such a way:
$\la_1^1=\la_1$, and for $1<j\le l$, $\la_1^j$ is defined to be
$$\align
&h_1\left(f_{l-j+2}^{(k_{l-j+2}^{j-1})}f_{l-j+1}^{(k_{l-j+1}^{j-1})}\cdots
f_1^{(k_1^{j-1})}\cdots f_l^{(k_l^1)}f_{l-1}^{(k_{l-1}^1)}\cdots f_1^{(k_1^1)}v\right)\\
=&\la_1^j\left(f_{l-j+2}^{(k_{l-i+2}^{j-1})}f_{l-j+1}^{(k_{l-i+1}^{j-1})}\cdots
f_1^{(k_1^{j-1})}\cdots f_l^{(k_l^1)}f_{l-1}^{(k_{l-1}^1)}\cdots f_1^{(k_1^1)}v\right)\\
=&\left(\la_1+\sum_{q=1}^{j-1}k^q_2-2\sum_{q=1}^{j-1}k^q_1\right)
\left(f_{l-j+2}^{(k_{l-j+2}^{j-1})}f_{l-j+1}^{(k_{l-j+1}^{j-1})}\cdots f_1^{(k_1^{j-1})}
\cdots f_l^{(k_l^1)}f_{l-1}^{(k_{l-1}^1)}\cdots f_1^{(k_1^1)}v\right);
\endalign$$
for $i>1$ and $j=1$, $\la_i^1$ is defined to be
$$
h_i\left(f_{i-1}^{(k_{i-1}^1)}\cdots f_{1}^{(k_{1}^1)}v\right)
=\la_i^1\left(f_{i-1}^{(k_{i-1}^1)}\cdots f_{1}^{(k_{1}^1)}v\right)
=(\la_i+k^1_{i-1})\left(f_{i-1}^{(k_{i-1}^1)}\cdots f_{1}^{(k_{1}^1)}v\right);
$$
for $i>1$ and $j>1$, $\la_i^j$ is defined to be
$$\align
&h_i\left(f_{i-1}^{(k_{i-1}^j)}\cdots f_{1}^{(k_{1}^j)}\cdots f_l^{(k_l^1)}\cdots
f_1^{(k_1^1)}v \right)=\la_i^j\left(f_{i-1}^{(k_{i-1}^j)}\cdots f_{1}^{(k_{1}^j)}
\cdots f_l^{(k_l^1)}\cdots f_1^{(k_1^1)}v\right)\\
=&\left(\la_i+\sum_{q=1}^{j}k^q_{i-1}+\sum_{q=1}^{j-1}k^q_{i+1}-2\sum_{q=1}^{j-1}k^q_{i}\right)
\left(f_{i-1}^{(k_{i-1}^j)}\cdots f_{1}^{(k_{1}^j)}\cdots f_l^{(k_l^1)}\cdots f_1^{(k_1^1)}v\right).
\endalign$$

Note that our definition is somewhat different from Littelmann's
definition in [2 \S 7].

Then we define the following two index set which are related to
$\la$ (comparing with Littelmann's definition of {\it S}$(\la)$ in
[2 \S 7]).
$$\Pi_{\la}=\Pi_{l, \la}:=\{K\in \Pi|~0\le k_i^j\le \la_i^j, \;
1\le i \le l, \; 1\le j \le l-i+1\}.$$
It is easy to see that $\Pi_{\la}$ is a finite set. We shall show
in Theorem 4.8 that the set $\{\theta^Kv|~K\in \Pi_{\la}\}$ forms
a basis of the $\bZ$-form of $V(\la)$.

For any $P\in \Pi'$, one has $P=(0, \cdots, 0, p_l$, $p_{l-1},
\cdots, p_1)$, if we set $p_0=0$, then we define
$$\Pi_{\la}':=\{P\in \Pi'|~p_{i}-p_{i-1}\le \la_i, \; 1\le i
<l,\},$$ and set $\la-\sum_{i=1}^lp_i\alpha_i=\la-P\alpha$ for
later use. We shall see in Theorem 4.7 that $\Pi_{\la}'$ is also a
finite set, and it becomes an index set of highest weights of
irreducible composition factors of $V(\la)$ to be viewed as a
$\g_{l-1})$-module.

{\bf 4.2.} Let $V$ be a $\U(\g_l)$-module. we say a vector $v\in V$ to
be a {\it primitive vector} of $V$, if there are two submodules $V_1,
V_2$ with $V_2\subset V_1\subseteq V$ such that $v\in V_1$, $v
\notin V_2$, and all $e_i$ with $1\le i\le l$ vanish the canonical
image of $v$ in $V_1/V_2$.

Let $V$ be a $\U(\g_l)$-module. According to [3], we can prove the
following lemma similarly.

{\bf 4.3.} {\bf {Lemma}}\;\;{\it {Let $w$ be a primitive vector of
weight $\la$ in $V$. Then $V$ has a composition factor isomorphic
to $V(\la)$.}}

Furthermore, one has the following lemma (cf. [1 \S 21.4.]).

{\bf 4.4.} {\bf {Lemma}}\;\;{\it {Let $\la =(\la_1, \la_2, \cdots,
\la_l)\in \Lambda^+$ be a dominant weight, and $v$ a maximal
vector of weight $\la$ of $V(\la)$. Then one has
$$f_i^{(\la_i +1)}v=0, \quad 1\le i \le l.$$}}

{\bf 4.5.} {\bf {Lemma}}\;\;{\it {Let $\la =(\la_1, \la_2, \cdots,
\la_l)\in \Lambda^+$ be a dominant weight. Let $V$ be a finite
dimensional $\U(\g_l)$-module generated by a maximal vector $v$ of
weight $\la$ of $V$. Then one has $V\simeq V(\la)$.}}

{\bf{Proof}\;\;} If $V$ is an irreducible $\U(\g_l)$-module. Then
$V\simeq V(\la)$. Otherwise, one has $$V=V(\la)\oplus M$$
according to the completely reducibility, because $V$ is a finite
dimensional $\U(\g_l)$-module. But $V$ is generated by a maximal vector,
it must be an indecomposable $\U(\g_l)$-module. This is a contradiction.

{\bf 4.6.} {\bf {Lemma}}\;\;{\it{Let $\la\! =\!(\la_1, \la_2,
\cdots, \la_l)\!\in\! \Lambda^+$ be a dominant weight, and
$P=(p_l, p_{l-1}$, $\cdots , p_1)$. Let $V(\la)$ be an irreducible
$\U(\g_l)$-module with maximal vector $v$. If there is an $i$ such
that $p_i-p_{i-1}>\la_i\ge 0$. Then one has
$$f_i^{(p_i)}f_{i-1}^{(p_{i-1})}\cdots f_1^{(p_1)}v=0.$$ }}

{\bf{Proof}\;\;}According to (2.2.7), one has
$$\align
&(f_i^{(p_i)}f_{i-1}^{(p_{i-1})})f_{i-2}^{(p_{i-2})}\cdots f_1^{(p_1)}v\\
=&\sum_{k=0}^{p_{i-1}}f_{i-1~i}^{(k)}f_{i-1}^{(p_{i-1}-k)}f_{i}^{(p_{i}-k)}
f_{i-2}^{(p_{i-2})}\cdots f_1^{(p_1)}v\\
=&\sum_{k=0}^{p_{i-1}}f_{i-1~i}^{(k)}f_{i-1}^{(p_{i-1}-k)}
f_{i-2}^{(p_{i-2})}\cdots f_1^{(p_1)}f_{i}^{(p_{i}-k)}v.
\endalign$$
Note that $k\le p_{i-1}$ and $0\le\la_i<p_i-p_{i-1}\le p_i-k$, the
above summation is zero by lemma 4.4.

Let $\la =(\la_1, \la_2, \cdots, \la_l)\in \Lambda^+$ be a dominant
weight. The finite-dimensional irreducible $\U(\g_l)$-module $V(\la)$
can be viewed as a $\U(\g_{l-1})$-module. It is no longer irreducible,
and can be decomposed into a direct sum of irreducible $\U(\g_{l-1})$
-module. The following theorem tell us how one can decompose it.

{\bf 4.7.} {\bf {Theorem}}\;\;{\it{Let $\la =(\la_1, \la_2, \cdots,
\la_l)\in \Lambda^+$ be a dominant weight. As a $\U(\g_{l-1})$-module,
the irreducible $\U(\g_l)$-module $V(\la)$ has the following direct sum
decomposition
$$V(\la)|_{\U(\g_{l-1})}=\bigoplus_{P\in \Pi_{\la}'}V\left(\left(
\la-P\alpha\right)_{\g_{l-1}}\right).$$}}

{\bf{Proof}\;\;}By definition, $\Pi_{\la}'$ is a finite set. Let
$|\Pi_{\la}'|=t$. We can arrange elements of $\Pi_{\la}'$
according to the ordering of $\Pi_{\la}'$ defined in \S 1.2. Then
one has $$\Pi_{\la}'=\{\text{\bf 0}=P_1 \prec P_2 \prec \cdots
\prec P_t\}.$$

Set $$M_{P_s}=\sum_{K\in \Pi, K \prec P_{s+1}}\bC \theta ^K v,
\quad 1\le s \le t-1,$$ where $v$ is a maximal vector of $V(\la)$
and $M_{P_t}=V(\la)$. Then one has
$$0\subseteq M_{P_1}\subseteq M_{P_2}\subseteq\cdots\subseteq
M_{P_t}=V(\la).$$

First of all, we show that $M_{P_s}$, $1\le s \le t$, is a
$\U(\g_{l-1})$-submodule of $V(\la)$. It does so when $s=t$. We
need only to consider cases of $1\le s <t$. For any $\theta^K v\in
M_{P_s}$ with $K\prec P_{s+1}$, it is still a weight vector, and
for any $h_i$ with $1\le i\le l$, one has by (2.1.3) $$h_i
\theta^K v=a_{i_K}\theta^K v \in M_{P_s}, \;\; \text {\rm
with}\;\; a_{i_K}\in \bZ.$$

By (3.4.1), $K=K_1+K_2$ with $K_1\in\Pi'$ and $K_2\in \Pi_{l-1}$.
Therefore, one has for any $f_i\in \U(\g_{l-1})$ with $1\le i \le
l-1$,
$$\align
f_i \theta^K v=&f_i \theta^{K_1+K_2} v=f_i
\left(\theta^{K_2}\theta^{K_1}\right) v\qquad\qquad\text{by
(3.4.2)}\\
=&\left(f_i \theta^{K_2}\right)\theta^{K_1} v=\left(
\sum_{K'\in \Pi_{l-1}}a_{_{K'}}\theta^{K'}\right)\theta^{K_1} v\\
=&\sum_{K'\in \Pi_{l-1}}a_{_{K'}}\theta^{K'+K_1} v, \;\; \text
{\rm with}\;\; a_{K'}\in \bZ.
\endalign$$
Note the fact that $K=K_1+K_2\prec P_{s+1}$, one has $K_1\prec
P_{s+1}$, and $K'+K_1\prec P_{s+1}$ for any $K'\in \Pi_{l-1}$.
Therefore,
$$f_i \theta^K v=\sum_{K'\in \Pi_{l-1}}a_{_{K'}}\theta^{K'+K_1}
v\in M_{P_s}.$$

Furthermore, one has for any $e_i$ with $1\le i \le l$,
$$\align
e_i \theta^K v=&e_i
f_1^{(k_1^l)}f_2^{(k_2^{l-1})}f_1^{(k_1^{l-1})}\cdots f_i^{
(k_i^{l-i+1})}f_{i-1}^{(k_{i-1}^{l-i+1})}\cdots
f_1^{(k_1^{l-i+1})} \cdots \\
&f_l^{(k_l^1)}f_{l-1}^{(k_{l-1}^1)}\cdots f_1^{(k_1^1)} v\\
=&\theta^Ke_iv+\sum_{n=1}^{l-i+1}f_1^{(k_1^l)}f_2^{(k_2^{l-1})}f_1^{(k_1^{l-1})}\cdots
f_i^{ (k_i^{n}-1)}(h_i-k_i^n+1)\\&f_{i-1}^{(k_{i-1}^{n})}\cdots
f_1^{(k_1^{n})} \cdots f_l^{(k_l^1)}f_{l-1}^{(k_{l-1}^1)}\cdots
f_1^{(k_1^1)}v\qquad\qquad\text{by
(2.1.5)}\\
=&\sum_{n=1}^{l-i+1}f_1^{(k_1^l)}f_2^{(k_2^{l-1})}f_1^{(k_1^{l-1})}\cdots
f_i^{ (k_i^{n}-1)}f_{i-1}^{(k_{i-1}^{n})}\cdots f_1^{(k_1^{n})}\cdots \\
&f_l^{(k_l^1)}f_{l-1}^{(k_{l-1}^1)}\cdots f_1^{(k_1^1)} a_n
v\qquad\qquad\qquad\qquad\qquad\quad\;\text{by
(2.1.3)}\\
=&\sum_{n=1}^{l-i+1}a_n\theta^{K-K_n}v,
\endalign$$
where
$$a_n=\la_i-k_i^n+1-2\sum_{d=1}^{n-1}k_i^d+\sum_{d=1}^{n}k_{i-1}^d+\sum_{d=1}^{n-1}k_{i+1}^d\in
\bZ,$$ and $K_n=(0, \cdots, 0, 1, 0, \cdots, 0)\in
\bN^{\frac{l(l+1)}{2}}$ with $1$ occurring in the place, where
$k_i^n$ lies in the corresponding $K$. Since $K-K_n\prec K\prec
P_{s+1}$, one has $$e_i \theta^K
v=\sum_{n=1}^{l-i+1}a_n\theta^{K-K_n}v\in M_{P_s}.$$ It shows that
$M_{P_s}$ is stable under actions of $e_i, h_i$ with $1\le i\le l$
and $f_i$ with $1\le i\le l-1$, and $M_{P_s}$ is a
$\U(\g_{l-1})$-module.

Secondly, we show that $\theta^{P_s}v$, $1\le s\le t$, are
primitive vectors in $V(\la)$ when it is viewed as a
$\U(\g_{l-1})$-module. Let $P_s=(0, \cdots 0, p_l,p_{l-1}, \cdots
, p_1)\in\Pi_{\la}'$. Then one has
$$\align
&e_1^{(p_1)}\cdots e_{l-1}^{(p_{l-1})} e_l^{(p_l)}\theta^{P_s}v\\
=&e_1^{(p_1)}\cdots e_{l-1}^{(p_{l-1})} e_l^{(p_l)}
f_l^{(p_l)}f_{l-1}^{(p_{l-1})}\cdots f_1^{(p_1)}v\\
=&e_1^{(p_1)}\cdots e_{l-1}^{(p_{l-1})}
\left(\sum_{k=0}^{p_l}f_l^{(p_l-k)}{{h_l-2p_l+2k}\choose
k}e_l^{(p_l-k)}\right)f_{l-1}^{(p_{l-1})}\cdots f_1^{(p_1)}v\\
=&e_1^{(p_1)}\cdots e_{l-1}^{(p_{l-1})}{{h_l}\choose p_l}
f_{l-1}^{(p_{l-1})}\cdots f_1^{(p_1)}v\\
=&e_1^{(p_1)}\cdots e_{l-1}^{(p_{l-1})}f_{l-1}^{(p_{l-1})}
{{h_l-p_{l-1}\al_{l-1}(h_l)}\choose p_l}f_{l-2}^{(p_{l-2})}
\cdots f_1^{(p_1)}v\qquad\text{by (2.1.4)}\\
=&e_1^{(p_1)}\cdots e_{l-1}^{(p_{l-1})}f_{l-1}^{(p_{l-1})}\cdots
f_1^{(p_1)}{{h_l-\sum_{k=1}^{l-1}p_k\al_k(h_l)}\choose
p_l}v\qquad\qquad\text{by (2.1.4)}\\
=&e_1^{(p_1)}\cdots e_{l-1}^{(p_{l-1})}f_{l-1}^{(p_{l-1})}\cdots
f_1^{(p_1)}{{\la_l+p_{l-1}}\choose p_l}v=\cdots=
\Pi_{k=1}^l{{\la_k+p_{k-1}}\choose p_k}v,
\endalign$$
where $p_0=0$, the second equality is by (2.1.2), and the last third
equality is because $\al_j(h_i)\ne 0$ if and only if $|i-j|\le 1$,
and $\al_j(h_{j\pm 1})=-1, \al_j(h_j)=2$ . Note that $p_{k}-p_{k-1}\le
\la_k$, one has $0\le p_k\le \la_k+p_{k-1}$, and ${{\la_k+p_{k-1}}\choose
p_k}\ne 0$ for all $1\le k\le l$, i.e. $e_1^{(p_1)}\cdots e_{l-1}^{(p_{l-1})}
e_l^{(p_l)}\theta^{P_s}v\ne 0$. This shows that $\theta^{P_s}v\ne 0$. By
our construction, it is easy to see that $\theta^{P_s}v\in M_{P_s}$ but
$\theta^{P_s}v\notin M_{P_{s-1}}$. Therefore, We need only to prove that
$e_i\theta^{P_s}v\in M_{P_{s-1}}$ for $1\le i \le l-1$, and then we can
conclude that $\theta^{P_s}v$ is a primitive vector in $V(\la)$. In fact,
one has for $1\le i\le l$
$$\align
e_i \theta^{P_s} v=&e_i f_l^{(p_l)}f_{l-1}^{(p_{l-1})}\cdots f_1^{(p_1)} v\\
=&\theta^{P_s}e_iv+f_l^{(p_l)}\cdots f_i^{(p_{i}-1)}(h_i-p_{i}+1)
f_{i-1}^{(p_{i-1})}\cdots f_1^{(p_1)}v \quad\quad\text{by (2.1.5)}\\
=&(\la_i-p_i+1+p_{i-1})f_l^{(p_l)}\cdots f_i^{(p_i-1)}f_{i-1}^{(p_{i-1})}
\cdots f_1^{(p_1)}v\qquad\quad\text{by (2.1.3)}\\
\endalign$$
Since $(0, \cdots, 0, p_l, \cdots,p_{i+1}, p_i-1, p_{i-1}, \cdots,
p_1)\prec P_{s}$, one has $e_i \theta^{P_s} v\in M_{P_{s-1}}$ as
required.

Thirdly, we show that $M_{P_s}=M_{P_{s-1}}+
\U(\g_{l-1})\theta^{P_s}v$. \lq\lq $\supseteq$\rq\rq is easy to be
proved by definition of $M_{P_s}$ and \S 3.4. Here we only prove
\lq\lq $\subseteq$\rq\rq. For any $K\in \Pi$ with $K\prec
P_{s+1}$, one has a unique decomposition $K=K_2+K_1$ with $K_2\in
\Pi'_{\la}$ and $K_1\in \Pi_{l-1}$. If $K\prec P_s$, then
$\theta^K v\in M_{P_{s-1}}$. Otherwise, when $P_s\preceq K\prec
P_{s+1}$, we must have $K_2=P_s$. Then
$$\theta^K v=\theta^{K_1+K_2}v=\theta^{K_1}\theta^{P_s}v\in \U(\g_{l-1})
\theta^{P_s}v$$
as required.

Finally, we show that $M_{P_s}/M_{P_{s-1}}\simeq V\left(\left(\la-
P_s\alpha\right)_{\g_{l-1}}\right)$. Let $w$ be the canonical image
of $\theta^{P_s}v$ in $M_{P_s}/M_{P_{s-1}}$. Then one has $M_{P_s}/
M_{P_{s-1}}\simeq \U(\g_{l-1})w$. Since $\theta^{P_s}v$ is a primitive
vector in $V(\la)$, $w$ becomes a maximal vector of weight $\left(\la-
P_s\alpha\right)_{\g_{l-1}}$. Note the fact that $V(\la)$ is a finite
dimensional module, and $M_{P_s}/M_{P_{s-1}}$ is also finite dimensional
and generated by a maximal vector $w$, we must have $M_{P_s}/M_{P_{s-1}}
\simeq V\left(\left(\la-P_s\alpha\right)_{\g_{l-1}}\right)$ by Lemma 4.5.

Using the complete reducibility, we complete the proof of Theorem 4.7.

The following theorem was proved in [2, Theorem 25].

{\bf 4.8.} {\bf {Theorem}}\;\;{\it{Let $v$ be a maximal vector of
$V(\la)$. Then $\{\theta^Kv|~K\in \Pi_{\la}\}$ forms a basis of
the $\bZ$-form of $V(\la)$. }}

{\bf{Proof}\;\;} We use induction on $l$. When $l=1$, one has for
any non-negative integer $m$ that $\{f_1^{(i)}v|~0\le i\le m\}$
forms a basis of the $\bZ$-form of $V(m)$ by Lemma 4.4. Assume
that our theorem holds for $l-1$, and then we have to show that
the theorem holds for $l$. Let us use the same notations as in the
proof of Theorem 4.7, and construct the bases of $M_{P_s}$ for
$1\le s \le t$. For $s=1$, one has $M_{P_1}\simeq V(\la_{\g_{l-1}})$
as $\U(\g_{l-1})$-module, and $\{\theta^Kv|~K\in \Pi_{{l-1}, {\la_
{\g_{l-1}}}}\}$ is a basis of $M_{P_1}$ by the induction hypothesis.
When $s=2$, note the following facts:

i) $\theta^{K+P_2}v\in M_{P_2}$ if $K\in \Pi_{{l-1},
{(\la-P_2\alpha)_{\g_{l-1}}}}$ by \S 3.4(3);

ii) the number of $\{\theta^K|K\!\in \Pi_{{l-1},
{(\la-P_2\alpha)_{\frak g_{l-1}}}}\}$ is equal to $\dim
V\left((\la-P_2\alpha)_{\frak g_{l-1}}\right)$ by the induction
hypothesis;

iii) $M_{P_2}/M_{P_1}\simeq V\left((\la-P_2\alpha)_{\frak
g_{l-1}}\right)$.
\newline
Therefore, we see that
$$\{\theta^Kv|~K\in \Pi_{{l-1},
{\la_{\frak g_{l-1}}}}\}\bigcup\{\theta^K \theta^{P_2}v=
\theta^{K+P_2}v |~K\in \Pi_{{l-1}, {(\la-P_2\alpha)_{\frak
g_{l-1}}}}\}$$ forms a basis of $M_{P_2}$.

In this way, the set of
$$\align
&\{\theta^Kv|~K\in \Pi_{{l-1}, {\la_{\frak
g_{l-1}}}}\}\bigcup\{\theta^{K+P_2}v
|~K\in \Pi_{{l-1}, {(\la-P_2\alpha)_{\frak g_{l-1}}}}\}\bigcup\\
&\qquad\qquad\quad\cdots\bigcup\{\theta^{K+P_t}v
|~K\in \Pi_{{l-1}, {(\la-P_t\alpha)_{\frak g_{l-1}}}}\}\\
\endalign$$
forms a basis of $M_{P_t}=V(\la)$. Note that elements in both the
above set and the set of $\{\theta^Kv|~K\in \Pi_{\la}\}$ are same,
this proves our theorem.

Denote by $\Pi(\la)$ the set of weights of the Weyl module
$V(\la)$. Let $P=(0, \cdots, 0, p_l$, $p_{l-1}$, $\cdots, p_1)\in
\Pi'_{\la}$. Then we say $P\alpha=\sum_ {i=1}^lp_i\alpha_i \ll
\sum_{i=1}^l a_i\alpha_i$ if and only if $p_l=a_l$ and $p_i\le
a_i$ for all $i=1, 2, \cdots, l-1$.

{\bf 4.9.} {\bf {Theorem}}\;\;{\it{Let $\mu \in \Pi(\lambda)$ be a
weight of $V(\la)$. Then the multiplicity $m_{\la}(\mu)$ of $\mu$
in $V(\la)$ is equal to}}
$$\align
m_{\la}(\mu)=\dim V(\la)_{\mu}&=\sum_{P\in\Pi_{\la}',\;P\alpha\ll
\la-\mu}\dim V\left(\left(\la-P\alpha\right)_{\g_{l-1}}\right)_{(\mu_{\g_{l-1}})}\\
&=\sum_{P\in \Pi_{\la}',\;P\alpha\ll \la-\mu}
m_{(\la-P\alpha)_{\g_{l-1}}}\left(\mu_{\g_{l-1}}\right).
\endalign$$

{\bf{Proof}\;\;} Let us use the same notations as in the proof of
Theorem 4.7, and let $\la-\mu=a_1\al_1+a_2\al_2+\cdots+a_l\al_l$
with all $a_i\ge 0$, $i=1, 2, \cdots, l$. Then the basis elements
of weight $\mu$ in $V(\la)$ are $\Cal M=\{\theta^Kv|~K=(k_1^l,
k_2^{l-1}, k_1^{l-1}, \cdots, k_i^{l-i+1}$, $k_{i-1}^{l-i+1}$,
$\cdots, k_1^{l-i+1}, \cdots, k_l^1, k_{l-1}^1$, $\cdots, k_1^1)
\in \Pi_{\la}$ with $k_l^1=a_l, k_{l-1}^1+k_{l-1}^2=a_{l-1},
\cdots, k_1^1+k_2^1+\cdots+k_l^1= a_1\}$, and the number of $\Cal
M$ is equal to $m_{\la}(\mu)$. If we divide $\Cal M$ into a disjoint
union of $\Cal M_i$, where $\Cal M_i=\{\theta^Kv|~K\in \Cal M$ with
$P_{i}\prec K\prec P_{i+1} \}$. From Theorem 4.8, we see that
$\Cal M_i\subseteq M_{P_i}$, and the number of $\Cal M_i$ is equal to
$m_{(\la-P_i\alpha)} \left(\mu_{\g_{l-1}}\right)$. Now Theorem 4.9
follows from Theorem 4.7.

\head {5. Examples}\endhead

{\bf 5.1.} When $l=2$, $\frak g_l$ is of type $A_2$. One has for
$\la=a\omega_1+b\omega_2=(a, b)\in \Lambda_+$ the following index
sets:
$$\align
&\Pi=\{(k^2_1, k^1_2, k^1_1)|~ k^1_1\le k^1_2\}\subseteq \bN^3,\\
&\Pi'=\{(0, k^1_2, k^1_1)|~ k^1_1\le k^1_2\}\subseteq \bN^3,\\
&\Pi_{\la}=\{(k^2_1, k^1_2, k^1_1)|~ k^1_1\le a, k_1^2\le
a+k_2^1-2k_1^1, k^1_2\le b+k_1^1\}\subseteq \Pi,\\
&\Pi_{\la}'=\{(0, p_2, p_1)|~ p_1\le a, p_2-p_1\le b\}\subseteq
\Pi'.\\
\endalign$$
In particular, if $\la=2\omega_1+3\omega_2=(2,3)$, then
$\Pi_{\la}=\{(k^2_1, k^1_2, k^1_1)\in \Pi|~ k^1_1\le 2, k_1^2\le
2+k_2^1-2k_1^1, k^1_2\le 3+k_1^1\}$, and $
\Pi_{\la}'=\{P_1=(0,0,0)\prec P_2=(0,1,0)\prec P_3=(0,2,0)\prec
P_4=$
\noindent $(0,3,0)\prec P_5=(0,1,1)\prec P_6=(0,2,1)\prec
P_7= (0,3,1)\prec P_8=(0,4,1)\prec P_9=(0,2,2)\prec
P_{10}=(0,3,2)\prec P_{11}=(0,4,2)\prec P_{12}=(0,5,2)\}$.

Moreover,
\newline
$M_{P_1}$ has basis $\{v, f_1v, f_1^{(2)}v\}$, and is isomorphic
to $V(2)$ as $\U(\g_1)$-modules;
\newline
$M_{P_2}/M_{P_1}$ has basis $\{f_2v, f_1f_2v, f_1^{(2)}f_2v,
f_1^{(3)}f_2v \}$, and is isomorphic to $V(3)$ as $\U(\g_1)$-modules;
\newline
$M_{P_3}/M_{P_2}$ has basis $\{f_2^{(2)}v, f_1f_2^{(2)}v,
f_1^{(2)}f_2^{(2)}v, f_1^{(3)}f_2^{(2)}v, f^{(4)}_1f_2^{(2)}\}$,
and is isomorphic to $V(4)$ as $\U(\g_1)$-modules;
\newline
$M_{P_4}/M_{P_3}$ has basis $\{f_2^{(3)}v, f_1f_2^{(3)}v,
f_1^{(2)}f_2^{(3)}v, f_1^{(3)}f_2^{(3)}v, f^{(4)}_1f_2^{(3)},
f^{(5)}_1f_2^{(3)}v\}$, and is isomorphic to $V(6)$ as $\U(\g_1)$-modules;
\newline
$M_{P_5}/M_{P_4}$ has basis $\{f_2f_1v, f_1f_2f_1v\}$, and is
isomorphic to $V(1)$ as $\U(\g_1)$-modules;
\newline
$M_{P_6}/M_{P_5}$ has basis $\{f_2^{(2)}f_1v, f_1f_2^{(2) }f_1v,
f_1^{(2)}f_2^{(2)}f_1v\}$, and is isomorphic to $V(2)$ as $\U(\g_1)$-modules;
\newline
$M_{P_7}/M_{P_6}$ has basis $\{f_2^{(3)}f_1v, f_1f_2^{(3) }f_1v,
f_1^{(2)}f_2^{(3)}f_1v, f_1^{(3)}f_2^{(3)}f_1v\}$, and is
isomorphic to $V(3)$ as $\U(\g_1)$-modules;
\newline
$M_{P_8}/M_{P_7}$ has basis $\{f_2^{(4)}f_1v, f_1f_2^{(4) }f_1v,
f_1^{(2)}f_2^{(4)}f_1v,f_1^{(3)}f_2^{(4)}f_1v,
f_1^{(4)}f_2^{(4)}f_1v\}$, and is isomorphic to $V(4)$ as $\U(\g_1)$-modules;
\newline
$M_{P_9}/M_{P_8}$ has basis $\{f_2^{(2)}f_1^{(2)}v\}$, and is
isomorphic to $V(0)$ as $\U(\g_1)$-modules;
\newline
$M_{P_{10}}/M_{P_9}$ has basis $\{f_2^{(3)}f_1^{(2)}v,\!
f_1f_2^{(3)}f_1^{(2)}v\!\}$, and is isomorphic to $V(1)$ as
$\U(\g_1\!)$-modules;
\newline
$M_{P_{11}}/M_{P_{10}}$ has basis $\{f_2^{(4)}f_1^{(2)}v,
f_1f_2^{(4)}f_1^{(2)}v, f_1^{(2)}f_2^{(4)}f_1^{(2)}v\}$, and is
isomorphic to $V(2)$ as $\U(\g_1)$-modules;
\newline
$M_{P_{12}}/M_{P_{11}}$ has basis $\{f_2^{(5)}f_1^{(2)}v,
f_1f_2^{(5)}f_1^{(2)}v, f_1^{(2)}f_2^{(5)}f_1^{(2)}v, f_1^{(3)}
f_2^{(5)}f_1^{(2)}v\}$, and is isomorphic to $V(3)$ as $\U(\g_1)$-modules.

Put all these elements together, we get a basis of
$V(2\omega_1+3\omega_2)$. Furthermore, one has
$V(2\omega_1+3\omega_2)|_{\U(\g_1)}\simeq\bigoplus_{i=1}^
{12}V(\la-P_i\alpha)|_{\U(\g_1)}$. It is known that
$m_{\la}(\mu)=3$ for $\mu=\omega_2$, and $\la-\mu=2\alpha_1
+2\alpha_2$. Using Theorem 4.9, one has
$m_{\la}(\mu)=m_{(\la-P_3\alpha)_{\g_1}} \left(\mu_{\g_1}\right)
+m_{(\la-P_6\alpha)_{\g_1}}\left(\mu_{\g_1}\right)+m_{(\la-P_9\alpha)_{\g_1}}
\left(\mu_{\g_1}\right)=m_4(0)+m_2(0)+m_0(0)=1+1+1=3$.

{\bf 5.2.} When $l=4$, $\frak g_l$ is of type $A_4$. One has for
$\la=a\omega_1+b\omega_2+c\omega_3+d\omega_4=(a, b, c, d)\in
\Lambda_+$ the following index sets:
$$\align
&\Pi=\{(k^4_1, k^3_2, k^3_1, k^2_3, k^2_2, k^2_1, k^1_4, k^1_3,
k^1_2, k^1_1)|~ k^3_1\le k^3_2, k^2_1\le k^2_2\le k^2_3, k^1_1
\le k^1_2\le k^1_3\le k^1_4 \}\\
&\;\;\quad\quad\subseteq \bN^{10},\\
&\Pi'=\{(0, \cdots, 0, k^1_4, k^1_3, k^1_2, k^1_1)|~ k^1_1\le
k^1_2\le k^1_3\le k^1_4 \}\subseteq \bN^{10},\\
&\Pi_{\la}=\{(k^4_1, k^3_2, k^3_1, k^2_3, k^2_2, k^2_1, k^1_4,
k^1_3, k^1_2, k^1_1)|~k^4_1\le a+k^3_2+k^2_2+k^1_2-2k^3_1-2k^2_1-2k^1_1,\\
&\qquad\quad k^3_2\le
b+k^3_1+k^2_1+k^1_1+k^2_3+k^1_3-2k^2_2-2k^1_2, k^3_1\le
a+k^2_2+k^1_2-2k^2_1-2k^1_1,\\
&\qquad\quad k^2_3\le c+k^2_2+k^1_2+k^1_4-2k^1_3, k^2_2\le
b+k^2_1+k^1_1+k^3_1-2k^1_2, k^2_1\le a+k^1_2-2k^1_1,\\
&\qquad\quad k^1_4\le d+k^1_3, k^1_3\le c+k^1_2, k^1_2\le b+k^1_1,
k^1_1\le
a\}\subseteq \Pi,\\
&\Pi_{\la}'=\{(0, \cdots, 0, p_4, p_3, p_2, p_1)|~ p_1\le a,
p_2-p_1\le b, p_3-p_2\le c, p_4-p_3\le d \}\subseteq
\Pi'.\\
\endalign$$
If we take $\la=\omega_1+\omega_2+\omega_3+\omega_4=(1,1,1,1)$,
then $\Pi_{\la}=\{(k^4_1, k^3_2, k^3_1, k^2_3, k^2_2, k^2_1,
k^1_4$, $ k^1_3, k^1_2, k^1_1)\in \Pi|~k^4_1\le
1+k^3_2+k^2_2+k^1_2-2k^3_1-2k^2_1-2k^1_1, k^3_2\le
1+k^3_1+k^2_1+k^1_1+k^2_3+k^1_3-2k^2_2-2k^1_2, k^3_1\le
1+k^2_2+k^1_2-2k^2_1-2k^1_1, k^2_3\le 1+k^2_2+k^1_2+k^1_4-2k^1_3,
k^2_2\le 1+k^2_1+k^1_1+k^3_1-2k^1_2, k^2_1\le 1+k^1_2-2k^1_1,
k^1_4\le 1+k^1_3, k^1_3\le 1+k^1_2, k^1_2\le 1+k^1_1, k^1_1\le 1
\}$, and $ \Pi_{\la}'=\{P_1=(0,\cdots, 0)\prec P_2=(0,\cdots,
0,1,0,0,0)\prec P_3=(0,\cdots,0,1,1,0,0)\prec
P_4=(0,\cdots,2,1,0,0)\prec P_5=(0,\cdots,0,1,1,1,0)\prec
P_6=(0,\cdots,0,2,1,1,0)\prec P_7= (0,\cdots,0,2,2,1,0)\prec
P_8=(0,\cdots,0,3,2,1,0)\prec P_9=(0,\cdots,0,1,1,1,1)\prec
P_{10}=(0,\cdots,0,2,1,1,1)\prec P_{11}=(0,\cdots,0,2,2,1,1)\prec
P_{12}=(0,\cdots,0,3,2,1,1)\prec P_{13}=(0,\cdots,0,2,2,2,1)\prec
P_{14}=(0,\cdots,0,3,2,2,1)\prec P_{15}=(0,\cdots,0,3,3,2,1)\prec
P_{16}=(0,\cdots,0,4,3,2,1)\}$.

Therefore, one has the following isomorphisms of  $\U(\g_3)$-modules:
$$\align
&M_{P_1}\simeq V(1,1,1),\qquad\;\;\quad M_{P_2}/M_{P_1}\simeq
V(1,1,2),\qquad
M_{P_3}/M_{P_2}\simeq V(1,2,0),\\
&M_{P_4}/M_{P_3}\simeq V(1,2,1),\;\quad M_{P_5}/M_{P_4}\simeq
V(2,0,1),\qquad M_{P_6}/M_{P_5}\simeq
V(2,0,2),\\
&M_{P_7}/M_{P_6}\simeq V(2,1,0),\;\quad M_{P_8}/M_{P_7}\simeq
V(2,1,1),\qquad M_{P_9}/M_{P_8}\simeq V(0,1,1),\\
&M_{P_{10}}/M_{P_9}\simeq V(0,1,2),\quad
M_{P_{11}}/M_{P_{10}}\simeq V(0,2,0),\;\quad M_{P_{12}}/M_{P_{11}}
\simeq V(0,2,1),\\
&M_{P_{13}}/M_{P_{12}}\simeq V(1,0,1),\quad
M_{P_{14}}/M_{P_{13}}\simeq V(1,0,2),\quad
M_{P_{15}}/M_{P_{14}}\simeq V(1,1,0),\\
&M_{P_{16}}/M_{P_{15}}\simeq V(1,1,1).\\
\endalign$$

Moreover, one has $V(\omega_1+\omega_2+\omega_3+\omega_4)|_
{\U(\g_3)}\simeq\bigoplus_{i=1}^ {16}V(\la-P_i\alpha)|_
{\U(\g_3)}$, and $m_{\la}(\mu)=8$ for $\mu=\omega_2+\omega_3=
(0,1,1,0)$ with $\la-\mu=\al_1+\alpha_2 +\alpha_3+\al_4$. Using
Theorem 4.9, one has
$$\align
m_{\la}(\mu)=&m_{(\la-P_2\alpha)_{\g_{3}}}
\left(\mu_{\g_{3}}\right)+m_{(\la-P_3\alpha)_{\g_{3}}}\left
(\mu_{\g_{3}}\right)+m_{(\la-P_5\alpha)_{\g_{3}}}
\left(\mu_{\g_{3}}\right)
+m_{(\la-P_9\alpha)_{\g_{3}}}\left(\mu_{\g_{3}}\right)\\
=&m_{(1,1,2)}(0,1,1)+m_{(1,2,0)}(0,1,1)+m_{(2,0,1)}(0,1,1)+m_{(0,1,1)}(0,1,1)\\
=&4+2+1+1=8.
\endalign$$

{\head {Acknowledgement}\endhead}

This work is supported in part by the National Natural Science
Foundation of China (10271088). The first named author is also
grateful to the Abdus Salam International Centre for Theoretical
Physics for its financial support and hospitability during his
visit.

\enddocument